\documentclass[11pt]{article}
\usepackage{amsmath, amsthm}
\usepackage{latexsym}
\usepackage{graphicx}
\makeatletter
\newfont{\footsc}{cmcsc10 at 8truept}
\newfont{\footbf}{cmbx10 at 8truept}
\newfont{\footrm}{cmr10 at 10truept}
\makeatother
\newtheorem{theorem}{\bf Theorem}

\newtheorem{lemma}{\bf Lemma}
\newtheorem{conjecture}{\bf Conjecture}

\title{Confirming Two Conjectures of Su and Wang}

\author{Yaming Yu\\
\small Department of Statistics\\[-0.8ex]
\small University of California\\[-0.8ex] 
\small Irvine, CA 92697, USA\\[-0.8ex]
\small \texttt{yamingy@uci.edu}}

\date{\small Mathematics Subject Classifications: 05A10, 05A20}

\begin{document}
\maketitle

\begin{abstract}
Two conjectures of Su and Wang (2008) concerning binomial coefficients are proved.  For $n\geq k\geq 0$ and $b>a>0$, we show 
that the finite sequence $C_j=\binom{n+ja}{k+jb}$ is a P\'{o}lya frequency sequence.  For $n\geq k\geq 0$ and $a>b>0$, we show 
that there exists an integer $m\geq 0$ such that the infinite sequence $\binom{n+ja}{k+jb},\ j=0, 1,\ldots$, is log-concave 
for $0\leq j\leq m$ and log-convex for $j\geq m$.  The proof of the first result exploits the connection between total 
positivity and planar networks, while that of the second uses a variation-diminishing property of the Laplace transform.
\end{abstract}

\section{Introduction}
A nonnegative sequence $u_i,\ i=0,1,\ldots$, is called {\it unimodal} if $u_0\leq\ldots\leq u_{m-1}\leq u_m\geq u_{m+1}\geq 
\ldots$ for some $m\geq 0$.  It is called {\it log-concave} (resp.\ {\it log-convex}), if $u_{i+1}^2\geq u_i u_{i+2}$ (resp.\ 
$u_{i+1}^2\leq u_i u_{i+2}$) for $i\geq 0$.  As is well-known, a log-concave sequence $u_i$ with no internal zeros (i.e., 
there exist no three indices $j<k<l$ such that $u_ju_l\neq 0$ but $u_k=0$) is unimodal.  Moreover, if a  polynomial $\sum_{i=0}^n u_i x^i$ with nonnegative coefficients has only real zeros, then the sequence $u_i,\ 0\leq i\leq n$, is log-concave with no internal zeros.  Unimodal, log-concave and log-convex sequences arise naturally in many problems in combinatorics and elsewhere; see \cite{B, LW} and \cite{Sa}--\cite{Y2}, for example.

Unimodality properties of sequences associated with Pascal's triangle have always been of interest (\cite{TZ, TZ2}).  
Recently, Su and Wang \cite{SW} have shown that the sequence of binomial coefficients located on a ray of Pascal's triangle 
is unimodal, as conjectured by Belbachir et al.\ \cite{BBS}.  At the end of \cite{SW}, the following new conjectures are 
proposed. 
\begin{conjecture}[\cite{SW}, Conjecture 2]
\label{conj0}
Let $n,\ k,\ a,\ b$ be integers such that $n\geq k\geq 0,\ b>a>0,$ and $k<b$. Define $C_j=\binom{n+ja}{k+jb},\ j=0, 1,\ldots$.  Then the polynomial $\sum_{j\geq 0} C_j x^j$ has only real zeros.
\end{conjecture}

\begin{conjecture}[\cite{SW}, Conjecture 3]
\label{conj}
Let $n,\ k,\ a,\ b$ be integers such that $n\geq k\geq 0$ and $a>b>0$.  Then there exists an integer $m\geq 0$ such that the 
sequence $\binom{n+ja}{k+jb},\ j=0, 1,\ldots$, is log-concave for $0\leq j\leq m$ and log-convex for $j\geq m$. 
\end{conjecture}
Note that in Conjecture \ref{conj0}, $C_j=0$ if $n+ja<k+jb$ by convention.  Also, in Conjecture \ref{conj}, $m=0$ is permitted, in which case the sequence $\binom{n+ja}{k+jb}$ is log-convex for all $j\geq 0$. 

In this work we confirm Conjectures \ref{conj0} and \ref{conj}.  Our proof of Conjecture \ref{conj0} (in Section 2) follows 
the combinatorial approach of Gessel and Viennot \cite{GV} and Brenti \cite{B2}.  In contrast, the proof of Conjecture 
\ref{conj} (in Section 3), which uses a variation-diminishing property of the Laplace transform, is analytic.  In the process 
of proving Conjecture \ref{conj}, we also obtain Theorem \ref{thm1}, which generalizes a result of Su and Wang (\cite{SW}, 
Proposition 1) that deals with the case $n=k=0$. 
\begin{theorem}
\label{thm1}
Assume $n\geq k\geq 0$ and $a>b>0$.  If $-1\leq k-(n+1)b/a\leq 0$ then the sequence $\binom{n+ja}{k+jb},\ j=0,1,\ldots,$ is 
log-convex.
\end{theorem}

As usual $\mathbf{N}$ denotes the set of positive integers and $\mathbf{Z}$ denotes the set of integers.

\section{Proof of Conjecture \ref{conj0}} 
Let us recall some useful terms.  An infinite matrix $W=(w_{ij})_{i,j\in \mathbf{N}}$ is called {\it totally 
positive} if every minor of $W$ is nonnegative.  A nonnegative sequence $u_i,\ i=0, 1,\ldots$, is 
called a {\it P\'{o}lya frequency sequence}, or PF sequence, if $(u_{j-i})_{i,j\in \mathbf{N}}$ ($u_i\equiv 0$ if $i<0$) is 
totally positive.  A finite sequence $u_i,\ i=0,\ldots, m,$ is a PF sequence if the infinite sequence $u_0, \ldots, u_m,\ 0, 
0, \ldots,$ is so.  The following connection between finite PF sequences and polynomials with real zeros is well-known; see 
Karlin \cite{K} for further notions and results concerning total positivity.

\begin{lemma}
\label{pf}
A nonnegative sequence $u_i,\ i=0,\ldots, n,$ is a PF sequence if and only if the polynomial $\sum_{i=0}^n u_i x^i$ has only real zeros.
\end{lemma}

Showing that a sequence is a PF sequence by definition can be nontrivial.  Nevertheless, it is possible to obtain remarkably 
simple proofs by exploiting the connection between total positivity and planar networks (\cite{B2, FZ, GV}).  In what follows, 
a {\it planar network} is a directed, acyclic, planar graph with no loops or multiple edges.  We allow the network to be 
infinite, but require that it is {\it locally finite}, i.e., there exist a finite number of paths between any two vertices.  
In addition, our network is associated with two sets of distinguished boundary vertices, one set on each side, and numbered 
from top to bottom as $s_i,\ i=1,2,\ldots$ (the sources) and $t_i,\ i=1,2,\ldots$ (the sinks) respectively.  (See \cite{B2} 
for a more precise formulation.)  Define the {\it path matrix} $W=(w_{ij})$ of such a network by 
$$w_{ij}={\rm\, the\ number\ of\ paths\ from\ } s_i\ {\rm to}\ t_j.$$
The following key lemma dates back to Karlin and McGregor \cite{KM} and Lindstr\"{o}m \cite{L}.  This and related techniques 
are used by Gessel and Viennot \cite{GV}, Stembridge \cite{Ste}, Sagan \cite{Sa} and Brenti \cite{B2} 
to tackle many combinatorial problems. 
\begin{lemma}
\label{lind}
The path matrix $W$ of a locally-finite planar network is totally positive.  Specifically, any $(I, J)$ minor of $W$ is equal to the number of families of vertex-disjoint paths that connect the sources labeled by $I$ with the sinks labeled by $J$.
\end{lemma}

We now construct a planar network with a particular path matrix suitable for applying Lemmas \ref{pf} and 
\ref{lind}.  Fix $n\geq k\geq 0$ and $b>a>0$.  Assume $k<b$ in addition, so that $C_j=\binom{n+aj}{k+bj}$ is indexed starting 
from $j=0$.  Let us specify the vertex set as 
$$V=\{(i,j):\ i, j\in \mathbf{Z},\ 0\leq i,\ 0\leq (b-a)i+bj\leq bn-ak\}.$$ 
For any two vertices $v_1=(i_1, j_1)$ and $v_2=(i_2, j_2)$ in $V$, we place an edge (oriented upwards and to the right) 
between $v_1$ and $v_2$ if $|i_1-i_2|+|j_1-j_2|=1$, i.e., the edge set is inherited from the square lattice $\mathbf{Z}\times 
\mathbf{Z}$.  Declare the vertices $s_i=(bi,\, (a-b)i),\ i=0,1,\ldots,$ as the sources, and $t_i=(k+bi,\, n-k+(a-b)i),\ 
i=0,1,\ldots$, as the sinks.  As an illustration, the special case $(n, k, a, b)=(4, 1, 1, 2)$ is displayed in Figure 1.

\begin{figure}
\begin{center}
\includegraphics[width=4.3in, height=3.3in]{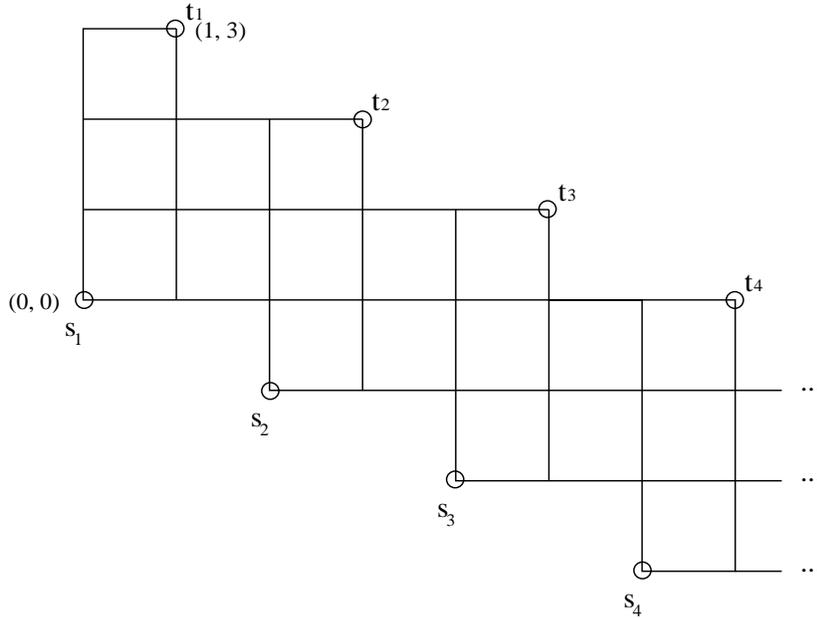}
\end{center}
\caption{The planar network corresponding to $(n, k, a, b)=(4, 1, 1, 2)$.}
\end{figure}

Evidently, for indices $i,\ j\geq 0$, if $j<i$ or $j>i+(n-k)/(b-a)$, then there are no paths from $s_i$ to $t_j$.  If $0\leq 
j-i\leq (n-k)/(b-a)$, then there are precisely $C_{j-i}=\binom{n+a(j-i)}{k+b(j-i)}$ such paths.  By Lemma \ref{lind}, the 
matrix $(C_{j-i})_{i,j\in \mathbf{N}}$ is totally positive; by Lemma \ref{pf}, Conjecture \ref{conj0} is valid. 

{\bf Remark.}  The Delannoy number (\cite{BS}) $D(n, k)$ counts the number of lattice paths from $(0, 0)$ to $(k, n)$ using 
only east, north and northeast steps; the recursion 
$$D(n, k)=D(n-1, k)+D(n, k-1)+D(n-1, k-1),\quad n, k\geq 1,$$
holds with the initial values $D(n, 0)=D(0, k)=1,\ n, k\geq 0$.  We have a result analogous to Conjecture \ref{conj0} for the 
Delannoy numbers.

\begin{theorem}
\label{delannoy}
Let $n,\ k,\ a,\ b$ be integers such that $n\geq k\geq 0$, $b>a>0$ and $k<b$.  Define $D_j=D(n-k+(a-b)j,\, k+bj)$.
Then the polynomial $\sum_{j\geq 0} D_j x^j$ has only real zeros.
\end{theorem}

Indeed, we can modify the planar network in the proof of Conjecture \ref{conj0} by adding all the edges from $(i, j)$ to 
$(i+1, j+1)$.  This new network then has $(D_{j-i})_{i,j\in \mathbf{N}}$ as its path matrix, 
and Theorem \ref{delannoy} follows from Lemmas \ref{lind} and \ref{pf} as before. 

It would be interesting to know whether results similar to Conjecture \ref{conj0} and Theorem \ref{delannoy} hold for the 
Stirling numbers of either kind, the Eulerian numbers, or their q-analogues.  Many matrices associated with these classical 
numbers are known or conjectured to be totally positive (Brenti \cite{B2}). 

\section{Proof of Conjecture \ref{conj}}
We shall analyze the quantity of interest as a Laplace transform.  A key tool is the following variation-diminishing property 
of the Laplace transform; see Karlin (\cite{K}, Chapter 5) for the precise statements and ramifications.
\begin{lemma}
\label{lem1}
Let $f(t)$ be a Borel-measurable function on $(0,\infty)$, and suppose the integral 
$$L(x)=\int_0^\infty f(t) e^{-xt}\, dt$$ 
converges absolutely for every $x\in (0, \infty)$.  Then the number of sign changes of $L(x)$ in $(0, \infty)$ is no more 
than the number of sign changes of $f(t)$ in $(0,\infty)$. 
\end{lemma}

Fix $n\geq k\geq 0,\ a>b>0$, and define 
$$g(x)=\log \frac{\Gamma(n+ax+1)}{\Gamma(k+bx+1)\Gamma(n-k+(a-b)x+1)},\quad x\geq 0,$$
where $\Gamma$ denotes Euler's gamma function.  Letting $\psi_1(x)=d^2 \log \Gamma(x)/dx^2$ as usual, and using the integral 
representation (\cite{AS}, p.\ 260)
\begin{equation}
\label{int}
\psi_1(x)=\int_0^\infty \frac{t e^{-x t}}{1-e^{-t}}\, dt,
\end{equation}
we get
\begin{align*}
g''(x) &=\int_0^\infty \frac{t}{1-e^{-t}} \left[a^2 e^{-(n+ax+1)t}-b^2 e^{-(k+bx+1)t} -(a-b)^2 e^{-(n-k+(a-b)x+1)t}\right]\, 
dt\\
       &=\int_0^\infty  a^2te^{-axt} \left[ \frac{e^{-(n+1)t}}{1-e^{-t}}- \frac{e^{-(k+1)ta/b}}{1-e^{-ta/b}} 
-\frac{e^{-(n-k+1)ta/(a-b)}}{1-e^{-ta/(a-b)}}\right]\, dt,
\end{align*}
where the second step uses two separate changes of variables.  For further simplification denote $u=k-(n+1)b/a,\ p=a/b,$ and 
$q=a/(a-b)$.  Note that $1/p+1/q=1$.  We obtain
\begin{equation}
\label{gx}
g''(x)=\int_0^\infty a^2 t e^{-axt-(n+1)t} h(t, u)\, dt
\end{equation}
with
$$h(t, u)=\frac{1}{1-e^{-t}}- \frac{e^{-(u+1)pt}}{1-e^{-pt}} -\frac{e^{uqt}}{1-e^{-qt}}.$$
It is easy to show that $\lim_{t\downarrow 0} h(t, u)=1/2$.  Also, $h(t, u)=O(e^{(n+1)t}),\ t\to\infty,$ for fixed $u$.  
By Watson's Lemma (see \cite{O}, for example), 
\begin{equation}
\label{watson}
g''(x) =\frac{a^2}{2(ax+n+1)^2} +o(x^{-2}),
\end{equation}
as $x\rightarrow\infty$.  This shows that $g(x)$ is asymptotically convex.  Note that a discrete version of this asymptotic 
convexity is obtained by Su and Wang (\cite{SW}, Theorem 1, part iii) using a different method.

Next, we examine the number of roots of $h(t, u)$ in $t\in (0,\infty)$ for fixed $u$.

\begin{lemma}
\label{lem2}
If $-1\leq u\leq 0$, then $h(t, u)>0$.
\end{lemma}
{\bf Proof.} It is easy to see that $h(t, u)$ is concave down in $u$.  Thus we only need to show $h(t, u)>0$ for $u=-1$ and 
$u=0$.  Let us assume $u=0$ since the case $u=-1$ can be obtained by switching the roles of $p$ and $q$.  We have
$$h(t, 0)=\frac{e^{-t}}{1-e^{-t}}- \frac{e^{-pt}}{1-e^{-pt}} -\frac{e^{-qt}}{1-e^{-qt}}.$$
Consider the function 
$$f(s)=\frac{s e^{-s}}{1-e^{-s}},\quad s>0.$$
It is easy to show that $f(s)$ strictly decreases in $s$.  Using this and $1/p+1/q=1$ we get
$$h(t, 0)=\frac{f(t)-f(tp)}{tp}+\frac{f(t)-f(tq)}{tq}>0. \qed $$
Note that Theorem \ref{thm1} follows directly from Lemma \ref{lem2} and expression (\ref{gx}).

\begin{lemma}
\label{lem3}
If $u\geq 0$ or $u\leq -1$ then $\partial h(t, u)/\partial t<0$.
\end{lemma}
{\bf Proof.} We may assume $u\geq 0$ since, as before, the case $u\leq -1$ can be obtained by switching the roles of $p$ and 
$q$.  We have
$$\frac{\partial h(t, 0)}{\partial t}=\frac{-e^{-t}}{(1-e^{-t})^2}+\frac{pe^{-pt}}{(1-e^{-pt})^2}+\frac{q 
e^{-qt}}{(1-e^{-qt})^2}.$$
It can be shown (details omitted) that the function 
$$l(s)=\frac{s^2 e^{-s}}{(1-e^{-s})^2},\quad s>0,$$  
strictly decreases in $s$.  Thus
\begin{equation}
\label{der1}
\frac{\partial h(t, 0)}{\partial t} =\frac{l(pt)-l(t)}{pt^2}+\frac{l(qt)-l(t)}{qt^2}<0.
\end{equation}

For $u>0$, it seems hard to determine the sign of $\partial h(t, u)/\partial t$ directly.  However, straightforward 
calculation gives  
$$\frac{\partial^2 h(t, u)}{\partial t\partial u} =\frac{p \left[ 1-pt-e^{-pt} -ptu(1-e^{-pt}) \right]}{e^{(u+1)pt} 
(1-e^{-pt})^2}+\frac{q \left[ (qt+1) e^{-qt} -1 -uqt(1-e^{-qt})\right]}{e^{-uqt} (1-e^{-qt})^2}.$$
In view of the simple inequalities $1-pt-e^{-pt}<0$ and $(qt+1) e^{-qt}-1<0$, we have 
\begin{equation}
\label{der2}
\frac{\partial^2 h(t, u)}{\partial t\partial u}<0,\quad u\geq 0.
\end{equation}
We obtain $\partial h(t, u)/\partial t<0,\ u\geq 0,$ from (\ref{der1}) and (\ref{der2}).  \qed

Lemmas \ref{lem2} and \ref{lem3} imply that, for fixed $u$, $h(t, u)$ has at most one root in $t\in (0, \infty)$.  Given 
(\ref{gx}) and Lemma \ref{lem1}, we know that $g''(x)$ changes sign at most once in $(0, \infty)$.  By (\ref{watson}), this 
possible change of sign is from negative to positive as $x$ increases from 0 to $\infty$.  Conjecture \ref{conj} is then 
proved when we restrict $x$ to be a non-negative integer. 

\section*{Acknowledgement}
The author would like to thank the Editor and an anonymous reviewer for their helpful comments.

\end{document}